\documentclass[a4paper,12pt,twoside]{article}
\usepackage{amsmath,amssymb,ifthen}
\usepackage[amsmath,thmmarks]{ntheorem}
\usepackage{paper}
\usepackage[all]{xy}
\usepackage{mathrsfs}
\usepackage{url}
\DeclareMathOperator{\PGL}{PGL}
\DeclareMathOperator{\rec}{rec}
\DeclareMathOperator{\Sp}{\mathbf{Sp}}
\DeclareMathOperator{\St}{\mathbf{St}}
\DeclareMathOperator{\LJ}{LJ}

\SelectTips{cm}{11}
\title{Note on weight-monodromy conjecture for $p$-adically uniformized varieties}
\author{Yoichi Mieda}

\def\headertitle{Weight-monodromy conjecture for $p$-adically uniformized varieties}

\begin{document}
\maketitle

\begin{firstfootnote}
 Graduate School of Mathematical Sciences, The University of Tokyo, 3--8--1 Komaba, Meguro-ku, Tokyo, 153--8914, Japan

 E-mail address: \texttt{mieda@ms.u-tokyo.ac.jp}

\end{firstfootnote}

\begin{abstract}
 We prove the weight-monodromy conjecture for varieties which are $p$-adically uniformized
 by a product of the Drinfeld upper half spaces. It is an easy consequence of Dat's work
 on the cohomology complex of the Drinfeld upper half space.
\end{abstract}

\section{Introduction}
Let $X$ be a proper smooth variety over a $p$-adic field $F$.
For a prime number $\ell\neq p$ and an integer $i$,
the absolute Galois group $\Gal(\overline{F}/F)$ acts on the $i$th $\ell$-adic \'etale cohomology
$H^i(X\otimes_F\overline{F},\overline{\Q}_\ell)$.
This action determines two filtrations on the cohomology; the weight filtration and
the monodromy filtration.
The weight-monodromy conjecture predicts that these two filtrations coincide up to shift by $i$.
This conjecture, due to Deligne \cite{MR0441965}, is widely open. It is known in the following cases:
\begin{enumerate}
 \item $X$ has good reduction over $\mathcal{O}_F$ (\cite{MR0340258}, \cite{MR601520}).
 \item $X$ is an abelian variety (\cite[Expos\'e IX]{SGA7}).
 \item $i\le 2$ (\cite{MR666636}, \cite{MR1423020}).
 \item $X$ is uniformized by the covering of the Drinfeld upper half space (\cite{MR2125735}, \cite{MR2224658}, \cite{MR2308851}). 
 \item $X$ is a set-theoretic complete intersection in a projective smooth toric variety (\cite{Scholze-perfectoid}).
\end{enumerate}
In this short note, we will slightly generalize the case (iv);
we will consider a variety $X$ which is uniformized by a product of the Drinfeld upper half spaces.

Our setting is as follows.
Let $F$, $F'$ be $p$-adic fields and $F''$ a $p$-adic field containing $F$ and $F'$.
Fix integers $d, d'\ge 1$ and put $G=\PGL_d(F)$, $G'=\PGL_d(F')$, respectively.
Let $\Omega_F=\Omega_F^{d-1}$ (resp.\ $\Omega_{F'}=\Omega_{F'}^{d'-1}$) denote the $d-1$-dimensional
(resp.\ $d'-1$-dimensional) Drinfeld upper half space.
To simplify the notation, we write $\Omega_F\times_{F''}\Omega_{F'}=(\Omega_F\otimes_F{F''})\times_{F''}(\Omega_{F'}\otimes_{F'}{F''})$.
For a discrete torsion-free cocompact subgroup $\Gamma\subset G\times G'$, the quotient
$\Omega_F\times_{F''}\Omega_{F'}/\Gamma$ becomes a projective smooth variety over $F''$.
Such a variety is said to be uniformized by $\Omega_F\times_{F''}\Omega_{F'}$.

The main theorem of this article is the following:

\begin{thm}\label{thm:main}
 Let $X$ be a projective smooth variety over $F''$ which is uniformized by
 $\Omega_F\times_{F''}\Omega_{F'}$. 
 Then, the weight-monodromy conjecture holds for $X$.
\end{thm}

Our strategy is the same as that in \cite{MR2224658}. First we determine the monodromy operator on
the cohomology complex
$R\Gamma_c((\Omega_F\times_{F''}\Omega_{F'})\otimes_{F''}\overline{F''},\overline{\Q}_\ell)$,
which is an object of the derived category of smooth $G\times G'$-representations.
Using this result, one can easily compute the cohomology of $X$,
from which the weight monodromy conjecture
is deduced.

Although Theorem \ref{thm:main} is stated for the product of two Drinfeld upper half spaces,
our argument also works for the product of more than two Drinfeld upper half spaces.
Further, as in \cite{MR2308851}, we may replace the Drinfeld upper half space
by its covering introduced by Drinfeld \cite{MR0422290}.
See Theorem \ref{thm:coh-covering}.

Interesting examples of varieties uniformized by a product of Drinfeld upper half spaces
are given by some unitary Shimura varieties (see \cite[Theorem 6.50]{MR1393439}).
By using our result, we can compute
the $\ell$-adic cohomology and the local Hasse-Weil zeta functions
of such Shimura varieties without any effort.
See the remark in the end of this note.

\medbreak
\noindent{\bfseries Acknowledgment}\quad
This work was supported by JSPS KAKENHI Grant Number 24740019.

\medbreak

\noindent\textbf{Notation}
Throughout this paper, $\ell$ denotes a prime number different from $p$.
Fix an isomorphism $\overline{\Q}_\ell\cong \C$ and identify them. 
All representations are considered over this field.

In the notation of $\ell$-adic \'etale cohomology, we omit the coefficient $\overline{\Q}_\ell$
and the base change to a separable closure. For example, 
$H^i(X\otimes_F\overline{F},\overline{\Q}_\ell)$ is written as $H^i(X)$.
As in Introduction, let $F$, $F'$ and $F''$ be $p$-adic fields with $F,F'\subset F''$.
For a rigid space $X$ (resp.\ $X'$) over $F$ (resp.\ $F'$), 
we write $X\times_{F''}X'$ for $(X\otimes_FF'')\otimes_{F''}(X'\otimes_{F'}F'')$.

\section{Proof}
For a subset $I$ of $\{1,\ldots,d-1\}$, an irreducible smooth representation $\pi_I$ of $G$
is naturally attached (see \cite[2.1.3]{MR2224658}). 
For example, $\pi_{\varnothing}$ is the Steinberg representation $\St_d$
and $\pi_{\{1,\ldots,d-1\}}$ the trivial representation $\mathbf{1}$.
For $0\le i\le d-1$, we write $\pi_{\le i}$ for $\pi_{\{1,\ldots,i\}}$.
Similarly, for $J\subset \{1,\ldots,d'-1\}$, consider an irreducible smooth representation $\pi'_J$
of $G'$.

\begin{lem}\label{lem:rep-theory}
 \begin{enumerate}
  \item For $I_1,I_2\subset \{1,\ldots,d-1\}$ and $J_1,J_2\subset\{1,\ldots,d'-1\}$, we have
	\[
	\Ext^i_{G\times G'}(\pi_{I_1}\boxtimes\pi'_{J_1},\pi_{I_2}\boxtimes\pi'_{J_2})
	=\begin{cases}\overline{\Q}_\ell& \text{if $i=\delta(I_1,I_2)+\delta(J_1,J_2)$,}\\
	 0& \text{otherwise.}\end{cases}
	\]
	Here $\delta(I_1,I_2)=\#(I_1\cup I_2)-\#(I_1\cap I_2)$.
  \item Let $I_1,I_2,I_3$ be subsets of $\{1,\ldots,d-1\}$ satisfying 
	$\delta(I_1,I_2)+\delta(I_2,I_3)=\delta(I_1,I_3)$.
	Take a non-zero element $\beta\in \Ext_G^{\delta(I_1,I_2)}(\pi_{I_1},\pi_{I_2})$.
	For $J_1,J_2\subset \{1,\ldots,d'-1\}$, the homomorphism
	\begin{align*}
	 &\Ext_{G\times G'}^{\delta(I_2,I_3)+\delta(J_1,J_2)}(\pi_{I_2}\boxtimes\pi'_{J_1},\pi_{I_3}\boxtimes\pi'_{J_2})\\
	 &\qquad\qquad \yrightarrow{(\beta\cup-)\boxtimes \id} \Ext_{G\times G'}^{\delta(I_1,I_3)+\delta(J_1,J_2)}(\pi_{I_1}\boxtimes\pi'_{J_1},\pi_{I_3}\boxtimes\pi'_{J_2})
	\end{align*}
	is an isomorphism.
  \item Let $\pi\boxtimes \pi'$ be an irreducible smooth representation of $G\times G'$.
	If it is not of the form $\pi_{I_0}\boxtimes \pi'_{J_0}$ with $I_0\subset \{1,\ldots,d-1\}$
	and $J_0\subset \{1,\ldots,d'-1\}$, then we have
	$\Ext^i_{G\times G'}(\pi_I\boxtimes\pi'_J,\pi\boxtimes\pi')=0$
	for every $I\subset \{1,\ldots,d-1\}$, $J\subset \{1,\ldots,d'-1\}$ and $i$.
 \end{enumerate}
\end{lem}

\begin{prf}
 For (i) and (ii), apply \cite[Th\'eor\`eme 1.3]{MR2224658} to the semisimple group $G\times G'$.
 The claim (iii) follows from \cite[Theorem 6.1]{MR1610824}, since the cuspidal supports of 
 $\pi_I\boxtimes \pi'_J$ and $\pi\boxtimes \pi'$ are different.
\end{prf}

Now we recall a result of Dat, which is crucial for our work.
In \cite{MR2224658}, he studied the cohomology complex 
$R\Gamma_c(\Omega_F)=R\Gamma_c(\Omega_F\otimes_F\overline{F},\overline{\Q}_\ell)$, which is
an object of the bounded derived category of smooth representations of $G$.
The Weil group $W_F$ of $F$ acts on $R\Gamma_c(\Omega_F)$.

\begin{thm}[\cite{MR2224658}]\label{thm:Dat}
 Fix a Frobenius lift $\varphi\in W_F$.
 \begin{enumerate}
  \item (\cite[Proposition 4.2.2]{MR2224658})
	There exists a unique isomorphism 
	\[
	 \alpha\colon R\Gamma_c(\Omega_F)\yrightarrow{\cong}\bigoplus_{i=0}^{d-1}\pi_{\le i}(-i)[-d+1-i]
	\]
	compatible with the actions of $\varphi$. It induces an isomorphism
	\[
	 \End\bigl(R\Gamma_c(\Omega_F)\bigr)\cong \bigoplus_{0\le i<j\le d-1}\Ext_G^{j-i}(\pi_{\le j},\pi_{\le i})(j-i).
	\]
  \item (\cite[Lemme 4.2.1]{MR2224658}) The monodromy operator $N\in \End(R\Gamma_c(\Omega_F))(-1)$
	on $R\Gamma_c(\Omega_F)$ is naturally determined.
	The image of $N$ under the isomorphism in (i) belongs to
	$\bigoplus_{0\le i\le d-2}\Ext_G^1(\pi_{\le i+1},\pi_{\le i})$.
	We denote it by $(\beta_i)_i$.
  \item (\cite[Proposition 4.2.7]{MR2224658}) For each $i$ with $0\le i\le d-2$, $\beta_i\neq 0$.
 \end{enumerate}
\end{thm}

The following theorem is an analogue of \cite[Th\'eor\`eme 1.1]{MR2224658}.

\begin{thm}\label{thm:coh}
 For subsets $I\subset \{1,\ldots,d-1\}$ and $J\subset \{1,\ldots,d'-1\}$, we have
 an isomorphism of Weil-Deligne representations of $F''$:
 \[
  \mathcal{H}^*\bigl(R\Hom\bigl(R\Gamma_c(\Omega_F\times_{F''}\Omega_{F'}),\pi_I\boxtimes \pi'_J\bigr)\bigr)\cong \rec_F(\pi_I)(\tfrac{d-1}{2})\vert_{W_{F''}}\otimes \rec_{F'}(\pi'_J)(\tfrac{d'-1}{2})\vert_{W_{F''}},
 \]
 where $\rec_F$ (resp.\ $\rec_{F'}$) denotes the local Langlands correspondence for $F$ (resp.\ $F'$).
 The functor $\mathcal{H}^*$ from $D^b(\overline{\Q}_\ell)$ to the category of 
 $\Z$-graded $\overline{\Q}_\ell$-vector spaces is given by $L^\bullet\mapsto \bigoplus_{i\in \Z}H^i(L^\bullet)$.

 If an irreducible smooth representation $\pi\boxtimes \pi'$ of $G\times G'$ is not of the form
 $\pi_I\boxtimes \pi'_J$, then
 $R\Hom(R\Gamma_c(\Omega_F\times_{F''}\Omega_{F'}),\pi\boxtimes \pi')=0$.
\end{thm}

\begin{prf}
 For simplicity, we only consider the cases $(I,J)=(\varnothing,\varnothing), (\varnothing,\{1,\ldots,d'-1\})$. Other cases can be treated similarly.

 First consider the case $(I,J)=(\varnothing,\varnothing)$.
 By Theorem \ref{thm:Dat} (i) and the K\"unneth formula, we have
 \[
  R\Gamma_c(\Omega_F\times_{F''}\Omega_{F'})\yrightarrow{\cong} \bigoplus_{i=0}^{d-1}\bigoplus_{j=0}^{d'-1}(\pi_{\le i}\boxtimes\pi'_{\le j})(-i-j)[-d-d'+2-i-j].
 \]
 By Lemma \ref{lem:rep-theory} (i), we have
 \[
  R\Hom\bigl(R\Gamma_c(\Omega_F\times_{F''}\Omega_{F'}),\pi_{\varnothing}\boxtimes \pi'_{\varnothing}\bigr)
 \cong \bigoplus_{i=0}^{d-1}\bigoplus_{j=0}^{d'-1}\Ext_{G\times G'}^{i+j}(\pi_{\le i}\boxtimes\pi'_{\le j},\pi_{\varnothing}\boxtimes \pi'_{\varnothing})(i+j),
 \]
 where $\Ext_{G\times G'}^{i+j}(\pi_{\le i}\boxtimes\pi'_{\le j},\pi_{\varnothing}\boxtimes \pi'_{\varnothing})$ is
 a one-dimensional vector space for each $i$ and $j$.
 Let $e_{0,0}\in \Hom_{G\times G'}(\pi_{\varnothing}\boxtimes \pi'_{\varnothing},\pi_{\varnothing}\boxtimes \pi'_{\varnothing})$ be the identity. Define $e_{i,j}\in \Ext_{G\times G'}^{i+j}(\pi_{\le i}\boxtimes\pi'_{\le j},\pi_{\varnothing}\boxtimes \pi'_{\varnothing})$ as the image of $e_{0,0}$ under the map
 \[
  (\beta_{i-1}\cup \cdots\cup \beta_0\cup -)\boxtimes (\beta'_{j-1}\cup \cdots\cup \beta'_0\cup -).
 \]
 Here, $(\beta'_j)\in \bigoplus_{0\le j\le d'-2}\Ext_{G'}^1(\pi'_{j+1},\pi'_j)$ denotes the image of
 $N\in \End(R\Gamma_c(\Omega_{F'}))(-1)$.
 By Lemma \ref{lem:rep-theory} (ii) and Theorem \ref{thm:Dat} (iii),
 $e_{i,j}$ is a basis of $\Ext_{G\times G'}^{i+j}(\pi_{\le i}\boxtimes\pi'_{\le j},\pi_{\varnothing}\boxtimes \pi'_{\varnothing})$. Now the monodromy operator on
 $R\Hom(R\Gamma_c(\Omega_F\times_{F''}\Omega_{F'}),\pi_{\varnothing}\boxtimes \pi'_{\varnothing})$ can be
 described explicitly:
 \[
  Ne_{i,j}=e_{i+1,j}+e_{i,j+1}.
 \]
 Therefore we conclude that
 \begin{align*}
  &\mathcal{H}^*\bigl(R\Hom(R\Gamma_c(\Omega_F\times_{F''}\Omega_{F'}),\pi_{\varnothing}\boxtimes \pi'_{\varnothing})\bigr)\cong \mathbf{Sp}_d(\tfrac{d-1}{2})\otimes \mathbf{Sp}_{d'}(\tfrac{d'-1}{2})\\
  &\qquad\qquad=\rec_F(\pi_\varnothing)(\tfrac{d-1}{2})\vert_{W_{F''}}\otimes \rec_{F'}(\pi'_\varnothing)(\tfrac{d'-1}{2})\vert_{W_{F''}}.
 \end{align*}
 We note that the right hand side concentrates in the degree $-d-d'+2$.

 Next assume that $I=\varnothing$ and $J=\{1,\ldots,d'-1\}$. Then, Lemma \ref{lem:rep-theory} (i)
 tells us that
  \begin{align*}
   &R\Hom\bigl(R\Gamma_c(\Omega_F\times_{F''}\Omega_{F'}),\pi_{\varnothing}\boxtimes \pi'_{\le d'-1}\bigr)\\
   &\qquad\cong \bigoplus_{i=0}^{d-1}\bigoplus_{j=0}^{d'-1}\Ext_{G\times G'}^{i+d'-1-j}(\pi_{\le i}\boxtimes\pi'_{\le j},\pi_{\varnothing}\boxtimes \pi'_{\le d'-1})(i+j).
 \end{align*}
 The $(i,j)$-component on the right hand side has degree $-d+1-2j$.
 Note that $\id\boxtimes (\beta'_j\cup-)$ induces the zero map on the right hand side.
 On the other hand, $(\beta_i\cup-)\boxtimes \id$ gives an isomorphism
 \[
  \Ext_{G\times G'}^{i+d'-1-j}(\pi_{\le i}\boxtimes\pi'_{\le j},\pi_{\varnothing}\boxtimes \pi'_{\le d'-1})\longrightarrow \Ext_{G\times G'}^{i+d'-j}(\pi_{\le i+1}\boxtimes\pi'_{\le j},\pi_{\varnothing}\boxtimes \pi'_{\le d'-1}).
 \]
 Therefore, $\mathcal{H}^*(R\Hom(R\Gamma_c(\Omega_F\times_{F''}\Omega_{F'}),\pi_{\varnothing}\boxtimes \pi'_{\le d'-1}))$ is isomorphic to
 \begin{align*}
  &\mathbf{Sp}_d(\tfrac{d-1}{2})\otimes \bigl(\overline{\Q}_\ell\oplus \overline{\Q}_\ell(1)\oplus\cdots \oplus \overline{\Q}_\ell(d'-1)\bigr)\\
  &\qquad\qquad=\rec_F(\pi_\varnothing)(\tfrac{d-1}{2})\vert_{W_{F''}}\otimes \rec_{F'}(\pi'_{\le d'-1})(\tfrac{d'-1}{2})\vert_{W_{F''}}.
 \end{align*}

 Finally, if $\pi\boxtimes \pi'$ is not of the form $\pi_I\boxtimes \pi'_J$, 
 \begin{align*}
  &R\Hom\bigl(R\Gamma_c(\Omega_F\times_{F''}\Omega_{F'}),\pi\boxtimes \pi'\bigr)\\
  &\qquad=\bigoplus_{i=0}^{d-1}\bigoplus_{j=0}^{d'-1}R\Hom(\pi_{\le i}\boxtimes\pi'_{\le j},\pi\boxtimes\pi')(i+j)[d+d'-2+i+j]=0
 \end{align*}
 by Lemma \ref{lem:rep-theory} (iii).
\end{prf}

\begin{cor}\label{cor:quotient-coh}
 Let $\Gamma$ be a discrete torsion-free cocompact subgroup of $G\times G'$.
 Let $m_{1,0}$ (resp.\ $m_{0,1}$, resp.\ $m_{1,1}$) be the multiplicity of
 ${\St_d}\boxtimes\mathbf{1}$ (resp.\ $\mathbf{1}\boxtimes \St_{d'}$, 
 resp.\ ${\St_d}\boxtimes {\St_{d'}}$) in the representation $C^\infty(G\times G'/\Gamma)$ of $G\times G'$.
 Then, we have a $W_{F''}$-isomorphism 
 \begin{align*}
  R\Gamma(\Omega_F\times_{F''}\Omega_{F'}/\Gamma)&\cong \Bigl(\bigoplus_{i=0}^{d-1}\bigoplus_{j=0}^{d'-1}\overline{\Q}_\ell(-i-j)[-2i-2j]\Bigr)\\
  &\qquad\qquad\oplus \biggl(\Sp_d(\tfrac{1-d}{2})[1-d]\otimes \Bigl(\bigoplus_{j=0}^{d'-1}\overline{\Q}_\ell(-j)[-2j]\Bigr)\biggr)^{m_{1,0}}\\
  &\qquad\qquad\oplus \biggl(\Bigl(\bigoplus_{i=0}^{d-1}\overline{\Q}_\ell(-i)[-2i]\Bigr)\otimes \Sp_{d'}(\tfrac{1-d'}{2})[1-d']\biggr)^{m_{0,1}}\\
  &\qquad\qquad\oplus (\Sp_d\otimes \Sp_{d'})^{m_{1,1}}(\tfrac{2-d-d'}{2})[2-d-d'].
 \end{align*}
 Moreover, the weight-monodromy conjecture holds for $\Omega_F\times_{F''}\Omega_{F'}/\Gamma$.
\end{cor}

\begin{prf}
 The proof is the same as \cite[Corollaire 4.5.1]{MR2224658}.
 By the fixed isomorphism $\overline{\Q}_\ell\cong\C$,
 we regard $C^\infty(G\times G'/\Gamma)$ as a representation over $\C$.
 Then, it is a unitary representation, and decomposes into the direct sum of irreducible smooth 
unitary representations of $G\times G'$ with finite multiplicities.
 Assume that $\pi_I\boxtimes \pi'_J$ appears in $C^\infty(G\times G'/\Gamma)$.
 Then it is unitary, and thus $\pi_I$ and $\pi'_J$ are unitary.
 Hence we conclude that $I$ (resp.\ $J$) is either $\varnothing$ or $\{1,\ldots,d-1\}$
 (resp.\ $\{1,\ldots,d'-1\}$). Note that the multiplicity of the trivial representation $\mathbf{1}\boxtimes \mathbf{1}$ in 
 $C^\infty(G\times G'/\Gamma)$ equals $1$.

 As in the proof of \cite[Corollaire 4.5.1]{MR2224658}, we have
 \begin{align*}
  &R\Gamma(\Omega_F\times_{F''}\Omega_{F'}/\Gamma)^\vee
  \cong \bigl(R\Gamma_c(\Omega_F\times_{F''}\Omega_{F'})\stackrel{\mathbb{L}}{\otimes}_{\overline{\Q}_\ell[\Gamma]}\overline{\Q}_\ell\bigr)^\vee\\
  &\qquad\cong R\Hom\bigl(R\Gamma_c(\Omega_F\times_{F''}\Omega_{F'}),C^\infty(G\times G'/\Gamma)\bigr)\\
  &\qquad=R\Hom\bigl(R\Gamma_c(\Omega_F\times_{F''}\Omega_{F'}),(\mathbf{1}\boxtimes \mathbf{1})\oplus ({\St_d}\boxtimes \mathbf{1})^{m_{1,0}}\\
  &\qquad\qquad\qquad\qquad\qquad\qquad\qquad\qquad\oplus (\mathbf{1}\boxtimes {\St_{d'}})^{m_{0,1}}\oplus ({\St_d}\boxtimes {\St_{d'}})^{m_{1,1}}\bigr).
 \end{align*}
 By using Theorem \ref{thm:coh} and taking dual, we obtain the desired description of
 $R\Gamma(\Omega_F\times_{F''}\Omega_{F'}/\Gamma)$.

 For the weight-monodromy conjecture, just note that 
 \[
 \Sp_d\otimes \Sp_{d'}\cong \bigoplus_{\substack{\lvert d-d'\rvert\le j\le d+d',\\j\equiv d+d'\bmod{2}}}\Sp_j
 \]
 (see \cite[(1.6.11.2)]{MR601520}).
\end{prf}

The argument above applies to the product of more than two Drinfeld upper half spaces
without any difficulty.
Furthermore, it is also valid even if we replace the Drinfeld upper half spaces by
its coverings introduced in \cite{MR0422290}. Namely, the following theorem holds.

\begin{thm}\label{thm:coh-covering}
 Let $D$ (resp.\ $D'$) be the central division algebra over $F$ (resp.\ $F'$) with invariant $1/d$
 (resp.\ $1/d'$). We denote $\mathcal{M}=\{\mathcal{M}_n\}$
 (resp.\ $\mathcal{M}'=\{\mathcal{M}'_{n'}\}$)
 the Drinfeld tower on which $D^\times$ (resp.\ $D'^\times$) acts.
 \begin{enumerate}
  \item Fix irreducible smooth representations $\rho$, $\rho'$ of $D^\times$, $D'^\times$,
	respectively. Let $\pi$ (resp.\ $\pi'$) be an irreducible smooth representation of $\GL_n(F)$
	(resp.\ $\GL_n(F')$) with the same central character as $\rho$ (resp.\ $\rho'$).
	\begin{enumerate}
	 \item If $\rho=\LJ_d(\pi)$ and $\rho'=\LJ_{d'}(\pi')$ (for the definition of $\LJ$, see \cite[\S 2]{MR2308851}),
	       then we have
	       \begin{align*}
	       &\mathcal{H}^*\bigl(R\Hom(R\Gamma_c(\mathcal{M}\times_{F''}\mathcal{M}')[\rho\boxtimes \rho']),\pi\boxtimes \pi'\bigr)\\
	       &\qquad\qquad\cong \rec_F(\pi)(\tfrac{d-1}{2})\vert_{W_{F''}}\oplus \rec_{F'}(\pi')(\tfrac{d'-1}{2})\vert_{W_{F''}}.
	       \end{align*}
	 \item Otherwise $R\Hom(R\Gamma_c(\mathcal{M}\times_{F''}\mathcal{M}')[\rho\boxtimes \rho']),\pi\boxtimes \pi')=0$.
	\end{enumerate}
	(See also \cite[Lemme 4.4.1]{MR2308851}.)
  \item Let $\Gamma$ be a discrete torsion-free cocompact subgroup of $\GL_d(F)\times \GL_{d'}(F')$,
	and $n, n'\ge 0$ integers.
	Then, $R\Gamma(\mathcal{M}_n\times_{F''}\mathcal{M}'_{n'}/\Gamma)$ can be computed as in
	\cite[p.~139--140]{MR2308851}. In particular, the weight-monodromy conjecture holds
	for $\mathcal{M}_n\times_{F''}\mathcal{M}'_{n'}/\Gamma$.
 \end{enumerate}
\end{thm}

\begin{prf}
 Use the result in \cite{MR2308851} in place of Theorem \ref{thm:Dat}.
\end{prf}

We may apply Theorem \ref{thm:coh-covering} to the unitary Shimura varieties appearing
in \cite[Theorem 6.50]{MR1393439}.
By the same method as in \cite[\S 3]{Shen-unitary}, one can compute the $\ell$-adic cohomology of
them using Theorem \ref{thm:coh-covering} (i).
This considerably simplifies the proof of the main result of \cite{Shen-unitary};
the study on test functions in \cite[\S 4--7]{Shen-unitary} is no longer needed.
The local Hasse-Weil zeta functions of such Shimura varieties can be computed directly.
The result is the same as in \cite[Corollary 7.4]{Shen-unitary} 
(but we do not need the assumption $r=1$).

\def\cftil#1{\ifmmode\setbox7\hbox{$\accent"5E#1$}\else
  \setbox7\hbox{\accent"5E#1}\penalty 10000\relax\fi\raise 1\ht7
  \hbox{\lower1.15ex\hbox to 1\wd7{\hss\accent"7E\hss}}\penalty 10000
  \hskip-1\wd7\penalty 10000\box7}
  \def\cftil#1{\ifmmode\setbox7\hbox{$\accent"5E#1$}\else
  \setbox7\hbox{\accent"5E#1}\penalty 10000\relax\fi\raise 1\ht7
  \hbox{\lower1.15ex\hbox to 1\wd7{\hss\accent"7E\hss}}\penalty 10000
  \hskip-1\wd7\penalty 10000\box7}
  \def\cftil#1{\ifmmode\setbox7\hbox{$\accent"5E#1$}\else
  \setbox7\hbox{\accent"5E#1}\penalty 10000\relax\fi\raise 1\ht7
  \hbox{\lower1.15ex\hbox to 1\wd7{\hss\accent"7E\hss}}\penalty 10000
  \hskip-1\wd7\penalty 10000\box7}
  \def\cftil#1{\ifmmode\setbox7\hbox{$\accent"5E#1$}\else
  \setbox7\hbox{\accent"5E#1}\penalty 10000\relax\fi\raise 1\ht7
  \hbox{\lower1.15ex\hbox to 1\wd7{\hss\accent"7E\hss}}\penalty 10000
  \hskip-1\wd7\penalty 10000\box7} \def\cprime{$'$} \def\cprime{$'$}
  \newcommand{\dummy}[1]{}
\providecommand{\bysame}{\leavevmode\hbox to3em{\hrulefill}\thinspace}
\providecommand{\MR}{\relax\ifhmode\unskip\space\fi MR }
\providecommand{\MRhref}[2]{%
  \href{http://www.ams.org/mathscinet-getitem?mr=#1}{#2}
}
\providecommand{\href}[2]{#2}

\end{document}